\documentclass[11pt]{amsart}
\usepackage{amsmath}
\usepackage{amssymb}

\newtheorem{theorem}{Theorem}[section]
\newtheorem{claim}[theorem]{Claim}
\newtheorem{lemma}[theorem]{Lemma}

\newtheorem{observation}[theorem]{Observation}
\newtheorem{corollary}[theorem]{Corollary}
\newtheorem{convention}[theorem]{Convention}

\theoremstyle{definition}
\newtheorem{definition}[theorem]{Definition}

\newtheorem{conjecture}[theorem]{Conjecture}

\theoremstyle{remark}
\newtheorem{remark}[theorem]{Remark}
\newtheorem{discussion}[theorem]{Discussion}

\newtheorem{conclusion}[theorem]{Conclusion}

\newcommand{\then}{{\underline{then}}}
\newcommand{\when}{{\underline{when}}}
\newcommand{\Then}{{\underline{Then}}}

\newcommand{\mn}{{\medskip\noindent}}
\newcommand{\sn}{{\smallskip\noindent}}

\newcount\skewfactor
\def\mathunderaccent#1#2 {\let\theaccent#1\skewfactor#2
\mathpalette\putaccentunder}
\def\putaccentunder#1#2{\oalign{$#1#2$\crcr\hidewidth
\vbox to.2ex{\hbox{$#1\skew\skewfactor\theaccent{}$}\vss}\hidewidth}}
\def\name{\mathunderaccent\tilde-3 }

\newcommand{\cA}{{\mathcal A}}

\newcommand{\cB}{{\mathcal B}}

\newcommand{\ga}{{\mathfrak a}}
\newcommand{\gI}{{\mathfrak I}}

\newcommand{\cF}{{\mathcal F}}

\newcommand{\bbP}{{\mathbb P}}
\newcommand{\cP}{{\mathcal P}}
\newcommand{\cQ}{{\mathcal Q}}

\newcommand{\bbZ}{{\mathbb Z}}

\newenvironment{PROOF}[2][\proofname.]
   {\begin{proof}[#1]\renewcommand{\qedsymbol}{\textsquare$_{\rm #2}$}}
   {\end{proof}}

\begin{document}

\title[Reflexive abelian groups]{Reflexive abelian groups and measurable
cardinals and full MAD families}

\author {Saharon Shelah}

\dedicatory{Dedicated to Walter Taylor}

\address{Einstein Institute of Mathematics\\
Edmond J. Safra Campus, Givat Ram\\
The Hebrew University of Jerusalem\\
Jerusalem, 91904, Israel\\
 and \\
 Department of Mathematics\\
 Hill Center - Busch Campus \\ 
 Rutgers, The State University of New Jersey \\
 110 Frelinghuysen Road \\
 Piscataway, NJ 08854-8019 USA}
\email{shelah@math.huji.ac.il}
\urladdr{http://shelah.logic.at}
\thanks{The author thanks Alice Leonhardt for the beautiful typing.
Research supported by German-Israeli Foundation for Scientific
Research and Development (Grant No. I-706-54.6/2001). 
Publication 904.}

\date{January 28, 2010}

\begin{abstract}
Answering problem (DG) of \cite{EM90}, \cite{EM02},
 we show that there is a reflexive group of cardinality which equals
 to the first measurable cardinal.  
\end{abstract}

\maketitle
\numberwithin{equation}{section}
\setcounter{section}{-1}

\centerline {Anotated Contents}
\bigskip

\noindent
\S0 Introduction
\bigskip

\noindent
\S1 A reflexive group above the first measurable
\begin{enumerate}
\item[${{}}$]   [We construct the abelian group $G_{\cA}$ for
every ${\cA} \subseteq {\cP}(\lambda)$ and find sufficient 
conditions for the existence of reflexive groups of cardinality at
least $\lambda$ among them.  From this we succeed to deduce
the existence of reflexive abelian groups of size the first measurable 
cardinality, answering a question from Eklof-Mekler book's.]
\end{enumerate}
\bigskip

\noindent
\S2 Arbitrarily large reflexive groups
\begin{enumerate}
\item[${{}}$]   [We show that it is ``very hard", not to have
reflexive groups of arbitrarily large cardinality.  E.g. after any set
forcing making the continuum above $\aleph_\omega$ not collapsing
$\aleph_\omega$, there are arbitrarily large reflexive groups.]
\end{enumerate}

\section {Introduction} 

\begin{definition}
\label{0.1} Let $G$ be an abelian group.
\begin{enumerate}
\item[(a)]  The dual of $G$ is the abelian group ${\rm Hom}(G,{\mathbb Z})$,
which we denote by $G^*$;
\item[(b)]  the double dual of $G$ is the abelian group ${\rm
Hom}(G^*,{\mathbb Z})$, which we denote by $G^{**}$.
\end{enumerate}
\end{definition}

There is the canonical homomorphism from $G$ into $G^{**}$, that is
$a \in G$ is
mapped to $F_a \in \text{ Hom}(G^*,\bbZ)$ defined by $F_a(f) = f(a)$.  The best
case, from our point of view, is when the canonical homomorphism is an
isomorphism.   There is a nice name for that phenomenon:

\begin{definition}
\label{0.2}
Let $G$ be an abelian group.  We say that $G$ is reflexive, if $G$ is
canonically isomorphic to $G^{**}$.
\end{definition}

Basic results about reflexive groups appear in Eklof and Mekler (see
\cite{EM90} and \cite{EM02} for a revised edition).  They present a
fundamental theorem of {\L}o\'s, generalized by Eda. {\L}o\'s theorem says
that $\lambda$ is smaller than the first
$\omega$-measurable\footnote{(set theorists 
call it the first measurable).} cardinal if
and only if the dual of the direct products of $\lambda$ copies of
$\mathbb Z$ is the direct sum of $\lambda$ copies of $\mathbb Z$.  
The inverse is always true.  It says that for all $\lambda$, the 
dual of the direct sum of
$\lambda$ copies of $\mathbb Z$ is the direct product of $\lambda$ copies of
$\mathbb Z$.  For $\lambda$ at least the first $\omega$-measurable,
{\L}o\'s's theorem just says the abelian group $\mathbb Z^\lambda$ is not
reflexive, Eda's theorem describes
${\rm Hom}({\Bbb Z}^\lambda,{\Bbb Z})$ in this case.  A direct consequence of
{\L}o\'s theorem is the existence of a lot of reflexive groups, but still
there is a cardinality limitation.  Let us describe the problem.  We use the
terminology of Eklof and Mekler.

\begin{definition}
\label{0.3}
Let $\mu$ be an infinite cardinal.
\begin{enumerate}
\item[(a)]  $\mu$ is measurable if there exists a non-principal
$\mu$-complete ultrafilter on $\mu$ and $\mu$ is uncountable
\item[(b)]  $\mu$ is $\omega$-measurable if there exists a
non-principal $\aleph_1$-complete ultrafilter on $\mu$.
\end{enumerate}
\end{definition}

We would like to clarify one important point.  Let $\mu$ be the first
$\omega$-measurable cardinal, and let $D$ be a non-principal
$\aleph_1$-complete ultrafilter on $\mu$.  It is well known that $D$ is also
$\mu$-complete.  So the first $\omega$-measurable cardinal is, in fact, the
first measurable cardinal.  It is easy to extend any non-principal
$\aleph_1$-complete ultrafilter on $\mu$ to an $\aleph_1$-complete
ultrafilter on $\lambda$.  So $\lambda$ is $\omega$-measurable for every
$\lambda \ge \mu$.

Let us summarize:

\begin{observation}
\label{0.4}
Let $\mu = \mu_{\text{first}}$ be the first measurable cardinal.
\begin{enumerate}
\item[(a)]  for every $\theta < \mu,\theta$ is not $\omega$-measurable
\item[(b)]  for every $\lambda \ge \mu,\lambda$ is $\omega$-measurable.
\end{enumerate}
\end{observation}

This terminology enables us to formulate the result that we need.  Recall
that ${\mathbb Z}^\theta$ is $\prod\limits_{i < \theta} {\mathbb Z}$ and
${\mathbb Z}^{(\theta)}$ is $\bigoplus\limits_{i<\theta} {\mathbb Z}$.  The
{\L}o\'s theorem deals with the existence of $\aleph_1$-complete
ultrafilters.  We will refer to the following corollary also known as {\L}o\'s
theorem:

\begin{corollary}
\label{0.5}
Let $\mu = \mu_{\text{first}}$ be the first measurable cardinal.
\begin{enumerate}
\item[(a)] for any $\theta < \mu$, ${\mathbb Z}^{(\theta)}$ is reflexive
(its dual being ${\mathbb Z}^\theta$).
\item[(b)] for every $\lambda \ge \mu$, ${\mathbb Z}^{(\lambda)}$ is not
reflexive. \hfill$\square_{\ref{0.5}}$
\end{enumerate}
\end{corollary}

The proof of \ref{0.5}(a) is based on the fact that every
$\aleph_1$-complete ultrafilter is principal.  So it does not work above
$\mu_{\rm first}$.  Naturally, we can ask - does there exist a reflexive
group of large cardinality, i.e., of cardinality $\ge \mu_{\rm first}$?
This is problem (DG) of Eklof-Mekler \cite{EM90}, \cite{EM02}.  We can
further ask

\begin{conjecture}
\label{ref.56.3}
There are reflexive abelian groups of arbitrarily large cardinalities.
\end{conjecture}

We thank Shimoni Garti, Daniel Herden and the referee for helpful
comments and corrections.

\section{A reflexive group above the first measurable cardinal}

We answer question (DG) of Eklof-Mekler \cite{EM90}.  There are reflexive
groups of cardinality not smaller than the first measurable.  Do we have it
for arbitrarily large $\lambda$, i.e. \ref{ref.56.3}?

This is very likely, in fact it follows (in ZFC) from $2^{\aleph_0} >
\aleph_\omega$ if some pcf conjecture holds.  See the next section.

\begin{theorem}
\label{ref.0A}
If $\mu = \mu_{\text{first}}$, the first measurable cardinal, \then\ there is
a reflexive $G \subseteq {}^\mu {\mathbb Z}$ of cardinality $\mu$.
\end{theorem}

\begin{PROOF}{\ref{ref.0A}}
By \ref{ref.21}(1) below (recall that $\mu$ is
measurable, so it is strong limit with cofinality greater than $\aleph_0$)
there are ${\mathcal A}_1,{\mathcal A}_2 \subseteq [\mu]^{\aleph_0}$
such that ${\mathcal A}_1 \subseteq {\mathcal A}^\perp_2$
and $\kappa^+({\mathcal A}^\perp_1) + \kappa^+({\mathcal A}^\perp_2)
\le \mu$, see Definition \ref{ref.0} below.  By claim \ref{ref.14}
below there is a $G$ as required.
\end{PROOF}

\begin{convention}  
\label{2.1A}
$\lambda \ge \aleph_0$ is fixed in this section (we need to fix $\lambda$ so
that ${\cA}^\perp$ is well defined).
\end{convention}

\begin{definition}
\label{ref.0}
\medskip

\noindent
\begin{enumerate}
\item For ${\mathcal A} \subseteq {\mathcal P}(\lambda)$, let
\begin{enumerate}
\item[$(a)$]  $\text{\rm id}({\mathcal A}) = \text{\rm id}_{\mathcal A}$
be the ideal of subsets of $\lambda$
generated by ${\mathcal A} \cup [\lambda]^{< \aleph_0}$;
\sn
\item[$(b$)]  ${\mathcal A}^\perp = \big\{u \subseteq \lambda:u \cap v$ finite
for every $v \in {\mathcal A}\big\}$;
\sn
\item[$(c)$]  $c\ell({\mathcal A}) = \big\{u \subseteq \lambda$: every
infinite $v \subseteq u$ contains some member of sb$({\mathcal A})$, 
see below$\big\}$;
\sn
\item[$(d)$]  {\rm sb}$({\mathcal A}) = \big\{u \subseteq \lambda:u$ is
infinite and is included in some member of ${\mathcal A}\ \big\}$.
\end{enumerate}
\item For ${\mathcal A} \subseteq {\mathcal P}(\lambda)$ let
\begin{enumerate}
\item[$(a)$]  $G_{\mathcal A}$ be the subgroup of $\mathbb Z^\lambda$
consisting of $\{f \in \mathbb Z^\lambda:{\rm supp}(f) \in 
\text{\rm id}({\mathcal A})\}$ where {\rm supp}$(f) = \{\alpha <
\lambda:f(\alpha) \ne 0_{\mathbb Z}\}$;
\sn
\item[(b)]  ${\mathbf j}_{\mathcal A}$ is the function from $G^*_{\mathcal A}
:= {\rm Hom}(G_{\mathcal A},{\mathbb Z})$ into  ${\mathbb Z}^\lambda$
defined by:
\[({\mathbf j}_{\mathcal A}(g))(\alpha) = g(\text{ch}_{\{\alpha\}})\]
where for $u\subseteq \lambda$, ch$_u = \text{ ch}_{\lambda,u}$ 
is the function with domain $\lambda$
mapping $\alpha$ to 1 if $\alpha \in u$ and to 0 if $\alpha \in \lambda
\setminus u$.
\end{enumerate}
\end{enumerate}
\end{definition}

\begin{definition}
\label{ref.0B}
\medskip

\noindent
\begin{enumerate}
\item[(a)] $\kappa^+({\mathcal A}) = \bigcup\{|u|^+:u \in {\mathcal A}\}$,
\sn
\item[(b)] $\kappa({\mathcal A}) = \bigcup\{|u|:u \in {\mathcal A}\}$.
\end{enumerate}
\end{definition}

\begin{claim}
\label{ref.7}
Let ${\mathcal A} \subseteq {\mathcal P}(\lambda)$.
\begin{enumerate}
\item
\begin{enumerate}
\item[(a)]  $\cA \subseteq c\ell({\mathcal A}) = c\ell(c\ell({\mathcal
A}))$,
\sn
\item[(b)]  ${\mathcal A}^\perp = \{u \subseteq \lambda:[u]^{\aleph_0} \cap
{\rm sb}({\mathcal A}) = \emptyset$, e.g. $u$ is finite$\}$
\sn
\item[(c)]  ${\mathcal A} \subseteq {\rm sb}({\mathcal A}) \cup [\lambda]^{<
\aleph_0}$
\sn
\item[(d)]  ${\mathcal A} \subseteq \text{\rm id}({\mathcal A})$
\end{enumerate}
\item
\begin{enumerate}
\item[(a)]  ${\mathcal A}^\perp = {\rm id}({\mathcal A}^\perp)$,
\sn
\item[(b)]  $({\mathcal A}^\perp)^\perp = c\ell({\mathcal A})$; note that
both include $[\lambda]^{< \aleph_0}$,
\sn
\item[(c)]  ${\mathcal A}^\perp = c\ell({\mathcal A}^\perp)$; note that
both include $[\lambda]^{< \aleph_0}$,
\sn
\item[(d)]  $\cA^\perp = (c\ell({\mathcal A}))^\perp$.
\end{enumerate}
\item If ${\mathcal A} \subseteq {\mathcal B} \subseteq {\mathcal
P}(\lambda)$, then ${\mathcal B}^\perp \subseteq {\mathcal A}^\perp$ and
$\kappa^+({\mathcal A}) \le \kappa^+({\mathcal B})$ and
$\kappa^+({\mathcal A})+\aleph_0=\kappa^+({\rm id}({\mathcal A}))$.
\end{enumerate}
\end{claim}

\begin{PROOF}{\ref{ref.7}}
1), 3) Obvious.

\noindent
2) {\bf Clause (a)}:\quad  Clearly ${\mathcal A}^\perp \subseteq
{\rm id}({\mathcal A}^\perp)$ by the definition of ${\rm id}$.  For the
other inclusion, as ${\mathcal A}^\perp$ includes all finite subsets of
$\lambda$ assume $u \in {\rm id}({\mathcal A}^\perp)$ is infinite hence
for some $n<\omega$ and infinite $u_0,\ldots,u_{n-1} \in {\mathcal A}^\perp$ we
have: $u \setminus \bigcup\{u_\ell:\ell < n\}$ is finite.  Hence

\[\begin{array}{l}
v \in {\mathcal A}\quad \Rightarrow\quad
(\forall \ell < n)(v \cap u_\ell\mbox{ is finite })\quad\Rightarrow\\\
(v \cap \cup\{u_\ell:\ell < n\}\mbox{ is finite })\quad\Rightarrow\quad
(v \cap u\mbox{ is finite })
\end{array}\]

\noindent
hence $u \in {\mathcal A}^\perp$.
\medskip

\noindent
{\bf Clause (b)}:\quad  Assume $u \in c \ell({\mathcal A})$
and $v \in {\mathcal A}^\perp$.  If $u \cap v$ is infinite then by
``$u \in c \ell({\mathcal A})$'' we know that $u \cap v$ includes some
member of ${\rm sb}({\mathcal A})$, but by ``$v \in {\mathcal A}^\perp$'' we
know that $u \cap v$ includes no member of sb$({\mathcal A})$,
contradiction.  So $u \cap v$ is finite.

Fixing $u \in c \ell({\mathcal A})$ and varying $v \in {\mathcal A}^\perp$
this tells us that $u \in (({\mathcal A})^\perp)^\perp$.  So we have
shown\footnote{recalling $[\lambda]^{< \aleph_0} \subseteq 
({\mathcal A}^\perp)^\perp$ is obvious!} $c \ell({\mathcal A}) 
\subseteq ({\mathcal A}^\perp)^\perp$.

Next if $u \subseteq \lambda$, $u\notin c \ell({\mathcal A})$ hence $u$ is
infinite then there is an infinite $v \subseteq u$ such that
$[v]^{\aleph_0} \cap {\rm sb}({\mathcal A}) =\emptyset$ hence $v$ is in
${\mathcal A}^\perp$, so $u$ includes an infinite member of ${\mathcal
  A}^\perp$ hence $u$ is not in $({\mathcal A}^\perp)^\perp$.  This shows
$u \notin c \ell({\mathcal A}) \Rightarrow u \notin ({\mathcal
A}^\perp)^\perp$.  So we get the desired equality.
\medskip

\noindent
{\bf Clause (c)}: \quad Similar to the proof of clause (b).
\medskip

\noindent
{\bf Clause (d)}:\quad Similar to the proof of clause (b).
\end{PROOF}

\begin{claim}
\label{ref.7B}
Let ${\mathcal A} \subseteq {\mathcal P}(\lambda)$.
\mn
\begin{enumerate}
\item If $\kappa^+({\mathcal A}) \le \mu_{{\rm first}} :=$ first
measurable cardinal, \then \, ${\mathbf j}_{\mathcal A}$ is an embedding
of $G^*_{\mathcal A}$ into ${\mathbb Z}^\lambda$ with its image
being $G_{\mathcal B}$ where ${\mathcal B} = {\mathcal A}^\perp$.
\sn
\item $G_{\mathcal A}$ is reflexive iff ${\rm id}({\mathcal A}) = c \ell({\mathcal
    A})$ and $\kappa^+({\mathcal A}) + \kappa^+({\mathcal A}^\perp) \le
  \mu_{{\rm first}}$.
\sn
\item $|G_{\mathcal A}| = \Sigma\{2^{|u|}:u \in {\rm id}({\mathcal A})\}
\in [\lambda,2^\lambda]$.
\sn
\item If $\kappa^+({\mathcal A})
  \le \mu$ then $\lambda \le |G_{\mathcal A}| = \lambda^{< \mu}$.
\end{enumerate}
\end{claim}

\begin{PROOF}{\ref{ref.7B}}
1) Clearly ${\mathbf j}_{\mathcal A}$ from
Definition \ref{ref.0}(2)(b) is a function from
  $G^*_{\mathcal A}$ into ${\mathbb Z}^\lambda$ and it is linear.  If $g \in
  G^*_{\mathcal A}$ and $u \in {\mathcal A}$ then by {\L}o\'s theorem (as
  $|u|<$ first measurable) necessarily $\{\alpha \in u:
g(\text{\rm ch}_{\{\alpha\}}) \ne 0\}$ is finite.  
So ${\rm supp}({\mathbf j}_{\mathcal A}(g)) \in {\mathcal
    A}^\perp$.  Together ${\mathbf j}_{\mathcal A}$ is a homomorphism from
  $G^*_{\mathcal A}$ into $G_{\mathcal B}$.  Also if ${\mathbf j}_{\mathcal
    A}(g_1) = {\mathbf j}_{\mathcal A}(g_2)$ but $g_1\ne g_2$ then for some
  $f \in G_{\mathcal A}$ we have $g_1(f) \ne g_2(f)$ and we can apply
  {\L}o\'s theorem for ${}^{{\rm supp}(f)}\mathbb Z$ to get a contradiction,
  hence $g_1 = g_2$ so we have deduced ``${\mathbf j}_{\mathcal A}$ is
  one-to-one''.  It is also easy to see that it is onto $G_{\mathcal B}$, so
  we are done. 

\noindent
2) First assume ${\rm id}({\mathcal A}) = c\ell({\mathcal A})$. By
  \ref{ref.7}(2)(c) we have ${\mathcal A}^\perp=c\ell({\mathcal A}^\perp)$.
  Applying part (1) to ${\mathcal A}$ and to ${\mathcal A}^\perp$ clearly if
  $\kappa^+({\mathcal A})+\kappa^+ ({\mathcal A}^\perp)\leq 
\mu_{\text{first}}$ we get $G_{\mathcal A}$ is canonically isomorphic to
  $(G^*_{\mathcal A})^*$.  Now $\bold j_{\cA}$ is an isomorphism from
  $G_{\cA}$ onto $G_{\cB} = G_{\cA^\perp}$ and $\bold j_{\cB}$ is an
  isomorphism from $G_{\cB}$ onto $G_{\cB^\perp}$, but by
  \ref{ref.7}(2)(c) by our assumption $\cB^\perp = (\cA^\perp)^\perp = c
  \ell(\cA) = \text{\rm id}(\cA)$.  So we have proved the ``if" implications.

If $\kappa^+({\mathcal A})>\mu_{\text{first}}$, then there
  is $u\in {\mathcal A}$ of cardinality $\ge \mu_{\text{first}}$, hence 
by \L os theorem we get $G_{\mathcal A}$ is not canonically
  isomorphic to $G^{**}_{\cA}$.

Lastly if ${\rm id}({\mathcal A}) \neq c\ell({\mathcal A})$ necessarily there
  is $u\in c\ell ({\mathcal A})\setminus {\rm id} ({\mathcal A})$ and let
$f = \text{ ch}_u$ so $u\in ({\mathcal A}^\perp)^\perp$ and $f$ defines
  a member of $(G_{{\mathcal A}^\perp})^*$ not
  ``coming from $G_{\mathcal A}$''.

\noindent
3),4)  Easy. 
\end{PROOF}

\begin{claim}
\label{ref.14}
A sufficient condition for the existence of a reflexive group $G$ of
cardinality $\ge \lambda$, in fact $\subseteq {\mathbb Z}^\lambda$ but
$\supseteq {\mathbb Z}^{(\lambda)}$ such that $|G| + |G^*| \le \lambda^{<
  \mu_{{\rm first}}}$, is $\circledast_{\lambda,\mu_{\text{\rm
first}}}$, when we define (for cardinals $\lambda \ge \mu$):
\medskip

\noindent
\begin{enumerate}
\item[$\circledast_{\lambda,\mu}$]  there are ${\mathcal A}_1,{\mathcal A}_2
  \subseteq [\lambda]^{\aleph_0}$ such that
\begin{enumerate}
\item[(a)]  ${\mathcal A}_1 \subseteq {\mathcal A}^\bot_2$, i.e.
\[u_1\in {\mathcal A}_1 \wedge u_2 \in {\mathcal A}_2\quad \Rightarrow\quad
u_1 \cap u_2\mbox{ is finite,}\]
\smallskip

\noindent
\item[(b)]  $\kappa^+({\mathcal A}^\bot_1) +\kappa^+({\mathcal A}^\perp_2)
\le \mu$.
\end{enumerate}
\end{enumerate}
\end{claim}

\begin{PROOF}{\ref{ref.14}}
  Let ${\mathcal A} = c \ell({\mathcal A}_1)$ and ${\mathcal B}
  = c\ell({\mathcal A}^\perp)$.  By \ref{ref.7}(2)(c) we have 
${\mathcal A}^\perp = {\mathcal B}$, 
and by \ref{ref.7}(2)(b), we have ${\mathcal B}^\perp =
  {\mathcal A}$ and lastly by \ref{ref.7}(1)(a) we have
  ${\mathcal A}_1 \subseteq {\mathcal A}$
  and by $\circledast_{\lambda,\mu}(a)$ we have
  ${\mathcal A}_2\subseteq {\mathcal A}_1^\perp$ but $\cA^\perp_1 = (c
\ell(\cA))^+ = \cA^\perp = c \ell(\cA^\perp) = \cB$ by the definitions
of $\cA$, by \ref{ref.7}(2)(d), hence by \ref{ref.7}(2)(c) and by the
definition of $\cB$; together we have ${\mathcal A}_2 \subseteq {\mathcal B}$.

Now ${\mathcal A}_1 \subseteq {\mathcal A}$ hence ${\mathcal A}^\perp
\subseteq {\mathcal A}^\perp_1$ so $\kappa^+({\mathcal A}^\perp) \le
\kappa^+({\mathcal A}^\perp_1)  \le \mu_{{\rm first}}$, also ${\mathcal A}_2
\subseteq {\mathcal B}$ hence ${\mathcal B}^\perp \subseteq {\mathcal
  A}^\perp_2$ hence $\kappa^+({\mathcal B}^\perp) \le \kappa^+({\mathcal
  A}^\perp_2) \le \mu_{{\rm first}}$.  But ${\mathcal A}^\perp = {\mathcal
  B}$ and  ${\mathcal B}^\perp  = {\mathcal A}$ and we have shown
$\kappa^+({\mathcal A}),\kappa^+({\mathcal B}) \le \mu_{{\rm first}}$.  So
by \ref{ref.7B}(1),(2)  $G_{\mathcal A},G_{\mathcal B}$ are reflexive and by
\ref{ref.7B}(4) the cardinality inequalities hold.
\end{PROOF}

\begin{claim}
\label{ref.21}
\medskip

\noindent
\begin{enumerate}
\item  If $\lambda > \aleph_0$ has uncountable cofinality, \then\ there are
  ${\mathcal A}_1,{\mathcal A}_2 \subseteq [\lambda]^{\aleph_0}$ such that
  ${\mathcal A}_1 \subseteq {\mathcal A}^\perp_2$ and $\kappa^+({\mathcal
    A}^\perp_1) + \kappa^+({\mathcal A}^\perp_2) \le \lambda$.
\smallskip

\noindent
\item  Assume $\lambda > \mu>\aleph_0$ 
and $S_1 \subseteq S^\lambda_{\aleph_0} =
\{\delta < \lambda:{\rm cf}(\delta) = \aleph_0\},S_2 = S^\lambda_{\aleph_0}
\setminus S_1$ are such that for every ordinal $\delta < \lambda$ of
cofinality $\ge \mu$, the set $\delta \cap S_\ell$ is stationary in
$\delta$ for $\ell=1,2$.  \Then\ for some ${\mathcal A}_1,{\mathcal A}_2
\subseteq [\lambda]^{\aleph_0}$, we have

\[\kappa^+({\mathcal A}^\perp_\ell) \le \mu\mbox{ for }\ell=1,2\quad \mbox{
  and }\quad {\mathcal A}_1 \subseteq {\mathcal A}^\perp_2.\]
\smallskip

\noindent
\item If ${\mathbf V} = {\mathbf L}$ (or, e.g. just $\neg \exists 0^\#$),
\then \, for every $\lambda > \aleph_0 = \mu$ the assumption of (2) holds.
\end{enumerate}
\end{claim}

\begin{PROOF}{\ref{ref.21}}  
1) For each $\delta \in S^\lambda_{\aleph_0}= \{\delta <
  \lambda:\delta$ limit of cofinality $\aleph_0\}$, let

\[{\mathcal P}_\delta=\{u\subseteq \delta:{\rm otp}(u) = \omega\mbox{ and
}\sup(u) =\delta\}.\]

\noindent
Let $S_1,S_2 \subseteq S^\lambda_{\aleph_0}$ be stationary disjoint subsets of
 $\lambda$ and let ${\mathcal A}_\ell = \bigcup\{{\mathcal P}_\delta:\delta
 \in S_\ell\}$ for $\ell=1,2$.  Now check.

\noindent
2) The same proof as the proof of part (1).

\noindent
3) Well known.
\end{PROOF}

\begin{remark}
\label{ref.24}
Also it is well known that
we can force an example as in \ref{ref.21}(2) for 
$\lambda= \mu_{{\rm first}}$, $\mu = \aleph_1$.

Without loss of generality ${\mathbf V}\models{\rm  GCH}$ and let $\theta =
{\rm  cf}(\theta) \le \mu_{{\rm first}},\theta > \aleph_0$.  Let
$\langle({\mathbb P}_\alpha,{{\name{\mathbb Q}}_\alpha}):\alpha \in {\rm
  Ord}\rangle$ be a full support iteration, ${\name{\mathbb Q}_\alpha}$ is
defined as follows: it is $\{f$: for some $\gamma < \aleph_\alpha$, $f \in
{}^\gamma\{1,2\}$, and for no increasing continuous sequence $\langle
\alpha_\varepsilon:\varepsilon < \theta\rangle$ of ordinals $< \gamma$ and
$\ell \in \{1,2\}$ do we have $\varepsilon < \theta \Rightarrow
f(\alpha_\varepsilon) = \ell\}$  if $\aleph_\alpha$ is regular, uncountable,
${\name{\mathbb Q}}_\alpha$ is trivial, $\{\emptyset\}$, otherwise.
\end{remark}

\begin{claim}
\label{ref.28}
Assume ${\mathbf V} = {\mathbf L}$ or much less: for every singular $\mu$
above $2^{\aleph_0}$ with countable cofinality, we have $\mu^{\aleph_0} =
\mu^+$ and $\square_\mu$.

\Then\ for every $\lambda$ there is a pair $({\mathcal A}_1,{\mathcal A}_2)$
as in $\circledast_{\lambda,\aleph_1}$ from \ref{ref.14}.
\end{claim}

\begin{PROOF}{\ref{ref.28}}
See Goldstern-Judah-Shelah \cite{GJSh:399}.
\end{PROOF}

\begin{remark}
\medskip

\noindent
\begin{enumerate}
\item  The assumption of Claim \ref{ref.28} holds in models with many
  measurable cardinals.
\smallskip

\noindent
\item  Note that if $\mu_1 \le \mu_2$ then clearly $\circledast_{\lambda,\mu_1}
\Rightarrow \circledast_{\lambda,\mu_2}$.
\end{enumerate}
\end{remark}

\section {Arbitrarily large reflexive groups}
In this section we shall show that it is ``hard'' to fail the
assumptions needed in the previous section in order to prove that there are
reflexive groups of arbitrarily large cardinality.  A typical result is
Conclusion \ref{ref.a45}.  Its proof uses parameters $\mathbf x$ (see
Definition \ref{ref.48}).  It is closed to an application in
\cite{Sh:460} to the Cantor discontinuum partition problem but
as the needed lemma \ref{ref.49} is only close to \cite{Sh:460}, we give a
complete proof in the appendix (the next section).
\bigskip

A characteristic conclusion is 

\begin{conclusion}
\label{ref.a45}
There is a reflexive subgroup $G$ of ${}^\lambda {\mathbb Z}$ if
$(*)_\mu$ below holds, moreover $G,G^*$ has cardinality $\in
[\lambda,\lambda^\mu]$, when
\medskip

\noindent
\begin{enumerate}
\item[$(*)_\mu$]  $\kappa$ is strong limit singular $<\mu_{{\rm first}}$ of
  cofinality $\aleph_0$ and $\kappa < \kappa^* < 2^\kappa$ and for
  no $\chi \ge 2^\kappa$ is there a subfamily ${\mathcal A} \subseteq
  [\chi]^{\kappa^*}$ of cardinality $> \chi$ the intersection of any two
  members is of cardinality $< \kappa$.
\end{enumerate}
\end{conclusion}

\begin{remark}
Alternatively, assume 
$\kappa = \aleph_0 < \kappa^*< 2^{\aleph_0}$, ${\mathfrak a}
  = 2^{\aleph_0}$.
\end{remark}

\begin{definition}
\label{ref.45}
\medskip

\noindent
\begin{enumerate}
\item We say that the triple $(\kappa,\kappa^*,\mu)$ is admissible
\when \, $\mu = \mu^\kappa$ (here usually $\mu = 2^\kappa$), 
$\kappa \le \kappa^* < \mu$ and the triple is
  $\lambda$-admissible for every $\lambda \ge \mu$, see below.
\smallskip

\noindent
\item The triple $(\kappa,\kappa^*,\mu)$ is $\lambda$-admissible \when \,
there is $\theta$ witnessing it which means:
\medskip

\noindent
\begin{enumerate}
\item[(a)] $\mu = \mu^\kappa,\kappa \le \kappa^* < \mu \le\lambda$,
\sn
\item[(b)] $\kappa^* \le \theta < \mu$,
\sn
\item[(c)] there is no family of more than $\lambda$ members of
$[\lambda]^{\ge \theta}$ such that the intersection of any two
has cardinality (strictly) less than $\kappa^*$.
\end{enumerate}
\item The triple $(\kappa,\kappa^*,\mu)$ is weakly $\lambda$-admissible
\when:
\mn
\begin{enumerate}
\item[(a)] as above, i.e. $\mu = \mu^\kappa,\kappa \le \kappa^* < \mu
\le\lambda$, 
\sn
\item[(b)] there is no family of more than $\lambda$ members of
$[\lambda]^\mu$ with any two of intersection of cardinality (strictly)
less than $\kappa^*$.
\end{enumerate}
\end{enumerate}
\end{definition}

\begin{remark}
\label{ref.45Y}
\medskip

\noindent
\begin{enumerate}
\item We may allow $(\kappa,\kappa^*)$ to be ordinals.
\item In the proof of \cite[3.8]{Sh:460}, ``$\theta$ witness
$(\kappa,\kappa^*,\mu)$ is $\lambda$-admissible'' was written
$\otimes^\theta_\lambda$.
\end{enumerate}
\end{remark}

\begin{claim}
\label{ref.46}
The triple $(\kappa,\kappa^*,\mu)$ is admissible \when \, at least one of the
following occurs:
\medskip

\noindent
\begin{enumerate}
\item[$(*)_1$]
\begin{enumerate}
\item[(a)] $\mu = 2^{\aleph_0} \geq \aleph_\delta > \kappa^* \ge \kappa =
  \aleph_0,\delta$ a limit ordinal
\item[(b)] for every $\lambda > \mu = 2^{\aleph_0}$ we 
have\footnote{recall
pp$_J(\theta) = \sup(\cup\{\text{pcf}_J(\bar\theta):\bar\theta =
\langle \theta_\varepsilon:\varepsilon \in S \rangle,
\theta_\varepsilon = \text{\rm cf}(\theta_\varepsilon) \in
(|S|,\theta)$ and $\theta = \lim_J\langle \theta_\varepsilon:
\varepsilon < \aleph_\alpha\rangle\})$ where $S = \text{ Dom}(J)$}

\begin{equation*}
\begin{array}{clcr}
\delta > \sup\{\alpha < \delta:&\text{ for some } \theta \in
(\mu,\lambda),\text{\rm cf}(\theta) = \aleph_\alpha \text{ and} \\
  &\text{\rm pp}_J(\theta) > \lambda 
\text{ for some } \aleph_\alpha \text{-complete ideal $J$ on } 
\aleph_\alpha\}
\end{array}
\end{equation*}
\end{enumerate}
\item[$(*)_2$] $\kappa > {\rm cf}(\kappa) = \aleph_0$ is strong
limit, $\delta$ a limit ordinal and we have:
\begin{enumerate}
\item[(a)]  $\mu = \mu^\kappa \ge \kappa^{+ \delta} > \kappa^* \ge \kappa$,
\sn
\item[(b)]  for every $\lambda > \mu$ we have

\begin{equation*}
\begin{array}{clcr}
\delta > \sup\{\alpha< \delta: &\text{ for some }\theta \in (\mu,\lambda),
\text{\rm cf}(\theta) = (\kappa^*)^{+ \alpha} \text{ and} \\
  &\text{\rm pp}_J(\theta) > \lambda \text{ for some }
(\kappa^*)^{+ \alpha} \text{-complete ideal on } (\kappa^*)^{+\alpha}\}.
\end{array}
\end{equation*} 
\end{enumerate}
\end{enumerate}
\end{claim}

\begin{remark}
\label{ref.47}  
In \ref{ref.46}, clause (b) of $(*)_2$ we can ask less
because in clause (c) of \ref{ref.45}(2) the intersection has cardinality
$< \kappa^*$ not just $< \theta$.
\end{remark}

\begin{PROOF}{\ref{ref.47}}
 Should be clear.  
\end{PROOF}

\begin{definition}
\label{ref.48}
1) The quintuple ${\mathbf x} = (X,c\ell,\kappa,\kappa^*,\mu)$ is 
{\em a parameter\/} \when:
\medskip

\noindent
\begin{enumerate}
\item[(a)] $c\ell:{\mathcal P}(X) \rightarrow {\mathcal P}(X)$,
\sn
\item[(b)] $\kappa \le \kappa^* \le \mu = \mu^\kappa$.
\end{enumerate}

\noindent
2)  The quintuple $\mathbf x$ is an admissible parameter \when \, in
addition:
\medskip

\noindent
\begin{enumerate}
\item[(c)] the triple $(\kappa,\kappa^*,\mu)$ is an admissible triple
(see Definition \ref{ref.45}(1) above).
\end{enumerate}

\noindent
3) We define

\[\begin{array}{ll}
{\mathcal P}^*_{\mathbf x} := \{A \subseteq X:&|A| = \mu \mbox{
and for every } B \subseteq A \mbox{ satisfying} \\
  &|B| = \kappa^* \mbox{ there is } B' \subseteq B,|B'| = \kappa
\mbox{ such that} \\
  &c \ell(B') \subseteq A, \mbox{ and } |c \ell(B')| = \mu\}
\end{array}\]

and

\[{\mathcal Q}^*_{\mathbf x} = \{B:B \subseteq X,|B| = \kappa \mbox{ and }
|c \ell(B)| = \mu\}\]

\noindent
and for $A \in {\mathcal P}^*_{\mathbf x}$ we define ${\mathcal
  Q}^*_{{\mathbf x}, A} = \{B \in {\mathcal Q}^*_{\mathbf x}:c \ell(B)
\subseteq A\}$.

\noindent
4) We say $\mathbf x$ is a strongly solvable parameter \when:
\medskip

\noindent
\begin{enumerate}
\item[(a),(b)]  as in part (1)
\sn
\item[(c)]  if $\bar h = \langle h^1_B,h^2_B:B \in {\mathcal Q}^*_{\mathbf
    x} \rangle$ and for every $B\in {\mathcal Q}^*_{\mathbf x}$ we have
  $h^\ell_B:c \ell(B) \rightarrow \mu$ for $\ell=1,2$ and $(\forall \alpha <
  \mu)(\exists^\mu \beta \in c \ell(B))(h^2_B(\beta) = \alpha)$, \then \,
 there is a function $h:X \rightarrow \mu$ such that:
\medskip

\noindent
\begin{enumerate}
\item[$\odot$]  if $A \in {\mathcal P}^*_{\mathbf x}$, so $|A|= \mu$ then
  for some $B \in {\mathcal Q}^*_{{\mathbf x},A}$ for every $\beta < \mu$
  the set $\{x \in c \ell(B):h^2_B(x) = \beta,h(x)=h^1_B(x)\}$ has
  cardinality $\mu$.
\end{enumerate}
\end{enumerate}
\medskip

\noindent
5) $\mathbf x$ is called solvable if above we restrict to the case
$h^2_B = h^1_B$.
\end{definition}

\begin{lemma}
\label{ref.49}
If ${\mathbf x} = (X,c \ell,\kappa,\kappa^*,\mu)$ is an admissible
parameter, \then\ $\mathbf x$ is strongly solvable.
\end{lemma}

\begin{PROOF}{\ref{ref.49}}
The proof is similar to \cite[3.8(2)]{Sh:460}, see a full proof
in the next section.
\end{PROOF}

We need the following for stating the main result:
\begin{definition}
\label{ref.53}
1) We say ${\mathcal A} \subseteq {\mathcal P}(\lambda)$ is
$(\sigma,\kappa^*,\mu)$-full in $\lambda$ \when \, ${\mathcal A} \subseteq
[\lambda]^\sigma$ and for every $A \in [\lambda]^{\kappa^*}$ 
 we have: $|A\cap B|\geq \sigma$ for at least
$\mu$ members $B$ of ${\mathcal A}$ \underline{or} $\sigma = \kappa^*$
and $\{B \in \cA:|B \cap A| \ge \sigma\}$ has cardinality $< \kappa^*$.

\noindent
2) We say ${\mathcal A}\subseteq [\lambda]^\sigma$ is
$(\sigma,\theta)$-MAD or $\theta$-MAD in $\lambda$ \when \,
$|{\mathcal A}|\geq \sigma$ and $B_1\neq B_2\in {\mathcal A}\Rightarrow
|B_1\cap B_2|< \theta$ and $B\in [\lambda]^\sigma \Rightarrow (\exists A\in
{\mathcal A}) (|A\cap B|\geq \theta)$.

\noindent
2A) If $\theta = \sigma$ we may omit $\theta$ writing ``MAD".  We may
omit ``in $\lambda$" and we may replace ``in $\lambda$" by ``in $A_*$".

\noindent
3) For $\theta\leq \sigma \leq \chi$ let 
${\mathfrak a}_{\chi,\sigma,\theta}={\rm Min}
\{|{\mathcal A}|:{\mathcal A}\subseteq [\chi]^\sigma$ is
$\theta$-MAD$\}$ and let ${\ga}_{\chi,\sigma} = {\ga}_{\chi,\sigma,\sigma}$.
\end{definition}

\begin{claim}
\label{ref.57}
1) Assume ${\mathcal A} \subseteq [\lambda]^\sigma$ is {\rm MAD},
i.e. $|\cA| \ge \sigma,A \ne B 
\in {\mathcal A} \Rightarrow |A \cap B| < \sigma$ and there is no $A \in
[\lambda]^\sigma$ such that $B \in {\mathcal A} \Rightarrow |B \cap A| <
\sigma$.  Then the family ${\mathcal A}$ is $(\sigma,\kappa^*,\mu)$-full 
(in $\lambda$) \when \, 
\mn
\begin{enumerate}
\item[$\boxplus_{\sigma,\kappa^*,\mu} $]  $\sigma \le 
\kappa^* < \mu$ and ${\ga}_{\kappa^*,\sigma} \ge \mu$.
\end{enumerate}
\mn
2) The statement $\boxplus_{\sigma,\kappa^*,\mu}$ holds when
at least one of the following occurs:
\mn
\begin{enumerate}
\item[$(*)_1$]   $\sigma = \aleph_0 \le \kappa^* < \mu = 2^{\aleph_0}$ and
  ${\mathfrak a} = 2^{\aleph_0}$ 
(or just ${\frak a}_{\kappa^*,\aleph_0}=2^{\aleph_0}$),
\sn
\item[$(*)_2$]  $\sigma$ is regular and for some strong limit singular
cardinal $\chi > \sigma$ of cofinality $\sigma$ we have $\chi \le \kappa^*
< \mu = 2^\chi$.
\end{enumerate}
\end{claim}

\begin{PROOF}{\ref{ref.57}}
1) Let $A \in [\lambda]^{\kappa^*}$, so if $\kappa^* >
\kappa$ then by ``$\cA$ is MAD" necessary $(\exists^{\ge \kappa^*}
B \in \cA)(B \cap A$ has cardinality $\ge \sigma$), hence
$(\exists^{\ge \kappa} B \in \cA)(B \cap A$ has cardinality $\ge \sigma$).  
Now ${\mathcal A}':= \{u\cap A:u \in {\mathcal A}$ and $u \cap
A$ has cardinality $\ge \sigma\}$ is a MAD
family of subsets of $A$ hence $|{\mathcal A}'|\geq
{\ga}_{\kappa^*,\sigma} \ge \mu$ as required.  
Note that $\ga_{\kappa^*,\sigma} \ge \ga_{\sigma,\sigma}$.

\noindent
2) \underline{Case 1}:  $(*)_1$ holds.

Obvious.
\medskip

\noindent
\underline{Case 2}:  We have $(*)_2$ so $\sigma,\chi,\kappa^*,\mu$ 
are as there.  Verifying $\boxplus_{\sigma,\kappa^*,\mu}$ the first
demand ``$\sigma \le \kappa^* < \mu$" is obvious - just check $(*)_2$,
but have to prove $\ga_{\kappa^*,\sigma} \ge \mu$; see Definition
\ref{ref.53}(3).  So assume $\cA \subseteq [\kappa^*]^\sigma$ is
$\sigma$-MAD in $\kappa^*$ and we should prove that $|\cA| \ge \mu$.

Let $A \in [\kappa^*]^\chi$ and
${\cA}' := \{u \cap A:u \in {\mathcal A}$ and 
$|u \cap A| = \sigma\}$ has cardinality $\ge\kappa^*$; clearly $\cA$ is
 a MAD subfamily of $[A]^\sigma$. But:
\medskip

\noindent
\begin{enumerate}
\item[$\odot_1$] there is a MAD family ${\mathcal A}_0 \subseteq [A]^\sigma$
of cardinality $\chi^\sigma = 2^\chi$,
\sn
\item[$\odot_2$]  if $u \in \cA'$ and even $u \in [A]^\sigma$ then
$|\{v \in {\mathcal A}_0:|v \cap u| \ge \sigma\}|\leq 2^\kappa$.
\end{enumerate}
\mn
Hence necessarily $|\cA'| = 2^\chi = \mu$,  hence $\{B \in \cA:A \cap
B$ of cardinality $\sigma\}$ has cardinality $\mu$ as required.
\end{PROOF}

We shall use the following definition for $\sigma = \aleph_0$ in 
the proof of the main result in
this section:

\begin{definition}
\label{ref.50}
Assume $\lambda$ is an infinite cardinal, ${\mathcal A} \subseteq
[\lambda]^\sigma$ a MAD family, $|{\mathcal A}| = \lambda^\sigma$
and $\bar u^* = \langle u^*_\alpha:\alpha < \lambda^\sigma\rangle$ enumerates
${\mathcal A}$ with no repetitions. 
For every $A \subseteq \lambda^\sigma$ we define set$(A) = 
\text{ set}(A,\bar u^*)$ as

\[
\bigcup\big\{u^*_\alpha:\mbox{ the set } u^*_\alpha
\cap(\bigcup\{u^*_\beta:\beta \in A\})\mbox{ is an infinite set}\big\} \cup
(\lambda \cap A).
\]
\end{definition}

\begin{claim}
\label{ref.35}
1) There is a reflexive group $G\subseteq {}^\lambda {\mathbb Z}$ of
  cardinality $\in [\lambda,\lambda^\mu]$ \when:
\medskip

\noindent
\begin{enumerate}
\item[(a)] $(\kappa,\kappa^*,\mu)$ is an admissible triple, $\mu <
\mu_{\text{\rm first}}$ first
\sn
\item[(b)] at least one of the following holds
\begin{enumerate}
\item[$(\alpha)$] ${\mathfrak a} = 2^{\aleph_0} = \mu$ and $\kappa = \aleph_0$
\sn
\item[$(\beta)$] $\kappa$ is strong limit singular of cofinality
$\aleph_0$ and $\mu = 2^\kappa$
\sn
\item[$(\gamma)$]   there is a MAD family
${\mathcal A} \subseteq [\mu]^{\aleph_0}$ which is
$(\aleph_0,\kappa^*,\mu)$-full, i.e. such that: if $A \in 
[\mu]^{\kappa^*}$ then

\[
|\{u \in {\mathcal A}:u \cap B\mbox{ is infinite }\}| = \mu.
\]
\end{enumerate}
\end{enumerate}

\noindent
2) Given $\mu$, for every $\lambda \ge \mu$ there are ${\mathcal A}_1,{\mathcal A}_2$ as in
$\circledast_{\lambda,\mu}$ of \ref{ref.14} provided that there are an
admissible triple, $(\kappa,\kappa^*,\mu)$ and a 
$(\aleph_0,\kappa^*,\mu)$-full MAD family
${\mathcal A}\subseteq [\lambda]^{\aleph_0}$.
\end{claim}

\begin{remark}  1) Concerning \ref{ref.35}(2) if $\kappa < \mu_{{\rm first}}$ 
then trivially \ref{ref.14} apply.

\noindent
2) Actually \cite[\S3]{Sh:668} deals essentially with equiconsistency
results for such properties.
\end{remark}

\begin{PROOF}{\ref{ref.35}} 
1) First there is a 
MAD family ${\mathcal A} \subseteq [\lambda]^{\aleph_0}$.  It is
$(\kappa,\kappa^*,\mu)$-full.  Why?  if assumption 
$(b)(\alpha)$ then by \ref{ref.57}
using $(*)_1$ in part (2) there; if $(b)(\beta)$ then 
by \ref{ref.57} using $(*)_2$ of part (2) there with
$\kappa$ here standing for $\chi$ there; of
course also (b)$(\gamma)$ implies this.
Second,  the result follows from part (2) and \ref{ref.14}.

\noindent
2) Without loss of generality $\lambda > \mu$, as otherwise the conclusion
is trivial. We use the Lemma \ref{ref.49}.

To apply it we shall choose $X,c \ell$ and let
${\mathbf x} = (X,c \ell,\kappa,\kappa^*,\mu)$ and show that the
demands there hold.  Let ${\mathcal A} \subseteq [\lambda]^{\aleph_0}$ be
a MAD family of cardinality $\lambda^{\aleph_0}$
which is $(\kappa,\kappa^*,\mu)$-full and without
loss of generality ${\mathcal A} \cap \lambda = \emptyset$, i.e. no 
$u \in \cA$ is a countable ordinal and let $\bar u^* = \langle
u^*_\alpha:\alpha < \lambda^{\aleph_0}\rangle$ list ${\mathcal A}$ with no
repetitions.

Recall, that by the claim's assumption 
$\mu = 2^\kappa$ and let $X =\lambda \cup {\mathcal A}$, i.e. if
$\alpha$ is an ordinal of cardinality $\aleph_0$ then $\alpha \notin \cA$.  We
define a function

\begin{equation*}
\begin{array}{clcr}
c\ell:&{\mathcal P}(\lambda\cup {\mathcal A}) \rightarrow
{\mathcal P}(\lambda \cup {\mathcal A}) \mbox{ by}: \\
  &c \ell(A) := A \cup\{B \in {\mathcal A}:B \cap {\rm set}(A,\bar u^*) 
\mbox{ is infinite}\}
\end{array}
\end{equation*}
\mn
where ${\rm set}(A,\bar u^*)$ is defined in Definition \ref{ref.50} with
$\aleph_0$ here standing for $\sigma$ there.

We shall prove
\begin{enumerate}
\item[$\otimes$]  the quintuple ${\mathbf x} =
  (X,c\ell,\kappa,\kappa^*,\mu)$ is an admissible parameter.
\end{enumerate}

We should check the demands of Definition \ref{ref.48}(1)(2),
\medskip

\noindent{\bf Clause (a)} there, i.e. $c \ell:{\mathcal P}(X)
\rightarrow {\mathcal P}(X)$ is trivial by our choices of $X,c \ell$
\medskip

\noindent{\bf Clause (b)} there says $\kappa \le \kappa^* \le \mu =
\mu^\kappa$, which is trivial.
\medskip

\noindent{\bf Clause (c)}: It says that $(\kappa,\kappa^*,\mu)$ is
admissible triple which holds by our assumption (a) of \ref{ref.35}.

So we can apply Lemma \ref{ref.49} hence
\begin{enumerate}
\item[$\oplus$]  $\mathbf x$ is strongly solvable, see Definition
  \ref{ref.48}(3).
\end{enumerate}
To apply it we should choose $\bar{h} = \langle h^1_u,h^2_u:u \in
{\mathcal Q}^*_{\mathbf x}\rangle$.
\smallskip

Given $u \in {\mathcal Q}^*_{\mathbf x}$, hence $u \in [X]^\kappa$ we let
$h^2_u:c \ell(u) \rightarrow \mu$ be such that
\[(\forall\alpha < \mu)(\exists^\mu \beta \in c \ell(u))[h^2_u(\beta) =
\alpha].\]
Let $h^1_u(x)$ be $h^2_u(x)$.

So by clause (c) of Definition \ref{ref.48}(4) there is a function
$h:X \rightarrow \mu$ satisfying $\odot$ from Definition
\ref{ref.48}(4). We define ${\mathcal A}_\ell := \{A \in {\mathcal A}:h(A)
= \ell\} \subseteq {\mathcal A}$ for $\ell=1,2$ and it suffices to check
that $({\mathcal A}_1,{\mathcal A}_2)$ are as required in
$\circledast_{\lambda,\mu}$ of Claim \ref{ref.14}.

First, as ${\mathcal A}_1,{\mathcal A}_2 \subseteq {\mathcal A}$ and
${\mathcal A}$ is $\subseteq [\lambda]^{\aleph_0}$ clearly ${\mathcal
  A}_1,{\mathcal A}_2 \subseteq [\lambda]^{\aleph_0}$.

Second, Clause (a) there says ``$u_1 \in {\mathcal A}_1 \wedge u_2 \in
{\mathcal A}_2 \Rightarrow u_1 \cap u_2$ finite'' and it holds as ${\mathcal
A}_1,{\mathcal A}_2$ are disjoint subsets of ${\mathcal A}$ which is a MAD
subset of $[\lambda]^{\aleph_0}$.

Third and lastly, clause (b) from $\circledast_{\lambda,\mu}$ of Definition
\ref{ref.14} says that $\kappa^+({\mathcal A}^\perp_\ell) \le \mu$.  So
towards a contradiction assume $A \subseteq \lambda,|A| = \mu$ and $A \in
{\mathcal A}^\perp_\ell$, i.e.

\[``u \in {\mathcal A}_\ell \Rightarrow A \cap u \text{ finite}"\]
Let

\[\cA' := \{u^*_\alpha \in {\mathcal A}:A \cap u^*_\alpha 
\text{ is infinite}\}.\]

Now if $A \in {\mathcal P}^*_{\mathbf x}$ then by the definition of ``$A \in
{\mathcal P}^*_{\mathbf x}$'' in \ref{ref.48}(3) there is $B \subseteq A$ which
belongs to ${\mathcal Q}^*_{\mathbf x,A}$ hence $|B| = \kappa$ and
recall $\kappa < \mu$.  Clearly $c \ell(B) \backslash
\lambda \subseteq \cA' \subseteq {\mathcal A}$ and by the choice of $h$ for
some such $B$ there is
$u^*_{\alpha_\ell} \in c \ell(B)$ satisfying
$h(u^*_{\alpha_\ell})= \ell$.  Also as
\[u^*_{\alpha_\ell} \in c\ell(B) \cap {\mathcal A} \subseteq c \ell(A) \cap
{\mathcal A},\]
clearly $h(u^*_{\alpha_\ell})=\ell$ and  
$u^*_{\alpha_\ell} \in {\mathcal A}_\ell$,
contradiction to ``$A\in {\mathcal A}_\ell^\perp$".  
So we are left with proving
\mn
\begin{enumerate}
\item[$\circledast$]  $A \in {\mathcal P}^*_{\mathbf x}$.
\end{enumerate}
This follows from the choice of ${\mathcal A} \subseteq
[\lambda]^{\aleph_0}$ as {\rm MAD} $(\kappa,\kappa^*,\mu)$-full.
\end{PROOF}

Note that it is very hard to fail $(\forall
\lambda)(\circledast_{\lambda,\mu_{{\rm first}}})$, e.g. easily

\begin{claim}
\label{ref.95}
1) If $\chi$ is strong limit (uncountable), $\mathbb P$ is a (set) forcing and,

\[\Vdash_{\mathbb P} ``2^{\aleph_0} > \chi \text{ and } \chi
\text{ is still a limit cardinal}"\]

where $\bbP$ has cardinality $\le \chi$ or at least 
satisfies the $\chi^+$-c.c., \then \, in $\mathbf V^{\mathbb P}$ the
triple $(\aleph_0,\chi,2^{\aleph_0})$ is admissible.

\noindent
1A) If $\bbP$ is a (set) forcing, $\kappa < \chi$ are strong limit
cardinals and $\Vdash_{\bbP} ``\chi$ is a limit cardinal $\kappa$ and is
strong limit cardinal of cofinality $\aleph_0$ and $\chi <
\kappa^{\aleph_0}$" and $\bbP$ satisfies the $\chi$-c.c., \then \, the
tuple $(\kappa,\chi,\kappa^{\aleph_0})$ is admissible. 

\noindent
2) If $\circledast_{\lambda,\mu}$ of \ref{ref.14} holds, $\mu =
\mu_{{\rm first}}$ or just $\mu$ is regular and $\mathbb P$ is a forcing
notion of cardinality $< \mu$, then we have $\circledast_{\lambda,\mu}$ in
$\mathbf V^{\mathbb P}$ also.
\end{claim}

\begin{PROOF}{\ref{ref.95}}
1) Without loss of generality there is $\delta$, a limit
ordinal such that $\Vdash_{\mathbb P} ``\mu = \aleph_\delta"$, and the
first demand of Definition \ref{ref.45}(1) and clause (a) of
\ref{ref.45}(2) hold.

By \cite{Sh:460}, or see \cite{Sh:829} in $\mathbf V$:
\mn
\begin{enumerate}
\item[$\odot_1$] for every $\lambda > \chi$ for some $\theta =
\theta_\lambda < \mu$ we have ${\rm cov}(\lambda,< \chi,<
\chi,\theta_\lambda) = \lambda$.
\end{enumerate}
\mn
This continues to hold in $\bold V^{\bbP}$ if we use
$\theta^1_\lambda = \theta_\lambda + (\text{cf}(\chi))^+$ if $\chi$ is
singular, $\theta'_\lambda = \theta_\lambda$ if $\chi$ is regular.

This is more than required in clauses (b),(c) of Definition \ref{ref.45}.
\medskip

\noindent
1A), 2) Easy.  
\end{PROOF}

\begin{remark}
\label{3.6A}
1) The holding of ``$\theta$ witness $(\kappa,\kappa^*,\mu)$ is
  $\lambda$-admissible'' is characterized in \cite[\S6]{Sh:410}.

\noindent
2) On earlier results concerning such problems and earlier
history see Hajnal-Juhasz-Shelah \cite{HJSh:249}.
\end{remark}

\section {Appendix: the proof of \protect{\ref{ref.49}}}

We are assuming ${\mathbf x} = (X,c\ell,\kappa,\kappa^*,\mu)$ is an
admissible parameter and we shall prove that it is strongly solvable.  In
Definition \ref{ref.48}(4) clauses (a),(b) hold trivially so it suffices to
prove clause (c). So let $\bar{h} = \langle h^1_B,h^2_B:B \in  
{\mathcal Q}^*_{\mathbf x}\rangle$ as there be given.

We prove by induction on $\lambda \in [\mu,|X|]$ that:
\medskip

\noindent
\begin{enumerate}
\item[$(*)_\lambda$] if $Z,Y$ are disjoint subsets of $X$ such that
$|Y| \le \lambda$, \then \, there are $h,Y^+$ such that
\begin{enumerate}
\item[(a)]  $Y \subseteq Y^+ \subseteq X \backslash Z$
\sn
\item[(b)]  $|Y^+| \le \lambda$
\sn
\item[(c)]  $h$ is a function from $Y^+$ to $\mu$
\sn
\item[(d)]  if $A \in {\mathcal P}^*_{\mathbf x}$, $\kappa^*\le\theta< \mu$,
the cardinal $\theta$ is a witness to $(\kappa,\kappa^*,\mu)$ being
  $\lambda$-admissible, $|A \cap Y^+| \ge \theta,|A \cap Z| < \mu$ and
$\beta < \mu$ \then \, 
$|\{x:h^2_B(x) = \beta \text{ and } h(x) = h^1_B(x)\}| = \mu$
for some $B \in \cQ^*_A$.
\end{enumerate}
\end{enumerate}
\medskip

\noindent
{\sc Case A:}\qquad  $\lambda = \mu$, so $|Y| \le \mu$.
\medskip

As $|Y| \le \mu = \mu^\kappa$, there is a set $Y^+$ of cardinality $\le
\mu$ such that $Y \subseteq Y^+ \subseteq X \backslash Z$ and
\mn
\begin{enumerate}
\item[$\odot_1$]  if $B \subseteq Y^+ \text{ and }
|B| \le \kappa \mbox{ and }|c\ell(B)| = \mu$ \then \, 
$c \ell(B) \backslash Z \subseteq Y^+$.
\end{enumerate}

Let

\[\begin{array}{ll}
{\mathcal P} =\{B \subseteq Y^+:&|B| \le \kappa\mbox{ and}\\
&(h^2_B)^{-1} (\{\beta\})\backslash Z\mbox{ has
cardinality $\mu$ for every }\beta < \mu\ \}.
\end{array}\]

\noindent
Clearly $|{\mathcal P}| \le |\{B:B \subseteq Y^+$ and $|B| \le \kappa\}|
\le |Y^+|^\kappa = \mu^\kappa = \mu$ and for every
$B \in {\mathcal P}$ and $\beta < \mu$ the set $(h^2_B)^{-1}(\{\beta\})
\backslash Z$ is included in $Y^+$ and has cardinality $\mu$.  So
$\langle(h^2_B)^{-1}(\{\beta\}) \backslash Z:
B \in {\mathcal P}$ and $\beta < \mu
\rangle$ is a sequence of 
$\mu$ subsets of $Y^+$ each of cardinality $\mu$.  Hence
there is a sequence $\langle C_{B,\beta}:B \in {\mathcal P},\beta <
\mu\rangle$ of pairwise disjoint sets such that $C_{B,\beta}
\subseteq (h^2_B)^{-1}(\{\beta\})$ and $|C_{B,\beta}| = \mu$.

Define a function $h$ from $Y^+$ to $\mu$ such that $h \restriction
C_{B,\beta} \subseteq h^1_B$ for $B \in {\mathcal P},\beta<\mu$ and

\[
h \restriction \big(Y^+ \backslash \bigcup\{C_{B,\beta}:B \in
{\mathcal P}\mbox{ and }\beta < \mu\}\big)\mbox{ is constantly zero.}
\]

Clearly clauses (a),(b),(c) of $(*)_\lambda$ holds.  For clause (d)
assume $A \in \cP^*_{\bold x}$ and $|A \cap Z| < \mu,\theta 
\in [\kappa^*,\mu)$ witness that
the tuple $(\kappa,\kappa^*,\mu)$ is $\mu$-admissible, and $|A
\cap Y^+| \ge \theta$.  Then by \ref{ref.48}(3) there is a set $B \in
\cQ^*_{\bold x,A}$, so $B \subseteq A,|B| \le \kappa$ and $|c \ell(B)|
=\mu$.  Clearly $B \in \cP$ and so clause (d) holds by the
choice of $h$.  So the function $h$ is as required.
\medskip

\noindent
{\sc Case B:}\qquad  $\lambda > \mu$.
\medskip

Let $\chi = \left( 2^\lambda \right)^+$ and choose $\langle N_i:i \le
\lambda \rangle$ an increasing continuous sequence of elementary submodels
of $({\mathcal H}(\chi),\in,<^*_\chi)$ such that
$X,c\ell,Y,Z,\lambda,\kappa,\kappa^*,\mu$ belong to the set $N_0, \mu + 1$
is included in $N_0$, the sequence $\langle N_i:i \le j \rangle$ belongs to
$N_{j+1}$ (when $j < \lambda$) and $\| N_i \| = \mu + |i|$.

Choose $\theta \in [\kappa^*,\mu)$ which witness that the triple
$(\kappa,\kappa^*,\mu)$ is $\lambda$-admissible.

We define by induction on $i < \lambda$, a set  $Y^+_i$ and a function
$h_i$ as follows:
\mn
\begin{enumerate}
\item[$\circledast$]  $(Y^+_i,h_i)$ is the $<^*_\chi$--first pair
  $(Y^\ast,h^\ast)$ such that:
\begin{enumerate}
\item[(a)] $Y^\ast\subseteq X \setminus (Z \cup\bigcup\limits_{j <
    i}Y^+_j)$
\sn
\item[(b)] $Y \cap N_i\backslash\bigcup\limits_{j < i} Y^+_j\backslash Z
  \subseteq X \cap N_i \backslash\bigcup\limits_{j < i} Y^+_j \backslash Z
  \subseteq Y^\ast$
\sn
\item[(c)] $|Y^\ast| \leq \mu + |i|$
\sn
\item[(d)] $h^\ast:Y^\ast \rightarrow \mu$
\sn
\item[(e)] $h^* \restriction ((h^2_B)^{-1} (\{\beta\})
\cap Y^*)$ coincides with $h^1_B$ on a set of cardinality $\mu$ 
for some $B\in {\mathcal Q}^*_{\bold x,A}$ and every $\beta < \mu$,
when for some $\theta'$:
\begin{enumerate}
\item[$(\alpha)$] $A \in {\mathcal P}^*_{\mathbf x}$,
\sn
\item[$(\beta)$]  $\kappa^\ast \le \theta' < \mu$, moreover $\theta'$ is a
  witness for the triple $(\kappa,\kappa^*,\mu)$ being $(\mu +
  |i|)$-admissible,
\sn
\item[$(\gamma)$] $|A \cap Y^*| \ge \theta'$,
\sn
\item[$(\delta)$] $|A \cap (Z \cup\bigcup\limits_{j < i} Y^+_j) | < \mu$.
\end{enumerate}
\end{enumerate}
\end{enumerate}

Note:  $(Y^+_i,h_i)$  exists by the induction hypothesis applied to the
cardinal $\lambda' := \mu + |i|$ and the sets $Z' := Z \cup
\bigcup\limits_{j < i} Y^+_j$ and $Y' := X \cap N_i \backslash
\bigcup\limits_{j < i}Y^+_j$ so we can carry out the induction.  
Also it is easy to prove by induction on $i$ that
\begin{enumerate}
\item[$\oplus$]  \begin{enumerate}
\item[(a)] $\langle (Y^+_j,h_j):j \le i\rangle \in N_{i+1}$
\sn
\item[(b)] $Y^+_j\subseteq N_{j+1}$
\end{enumerate}
\end{enumerate}
\mn
[Why?  First we show $\langle (Y^+_j,h_j):j<i\rangle \in N_{i+1}$ as
the induction can be carried inside $N_{i+1}$.  Now $Y^+_i,h_i \in N_{i+1}$
 as we always have chosen ``the $<^*_\chi$-first'', so clause (a)
above holds.  As for $Y^+_i \subseteq N_{i+1}$; i.e. 
clause (b) note that $|Y^+_j|=\mu+|i|$ and
$(\mu+1)\subseteq N_{j+1}, j+1\subseteq N_{j+1}$ by the choice of
$\langle N_i:i<\lambda\rangle$.]  

Let $Y^+ = \bigcup\limits_{i < \lambda} Y^+_i$ and $h = \bigcup\limits_{i <
\lambda} h_i$. Clearly $Y \subseteq N_\lambda 
= \bigcup\limits_{i < \lambda} N_i$
 as $Y\in N_0,i < \lambda \Rightarrow i \subseteq N_i \subseteq N$ 
and $|Y|=\lambda$ so $\lambda\in N_0$ and 
$\lambda\subseteq N_\lambda$ ,
hence by requirement (b) of $\circledast$ clearly $Y \subseteq Y^+$, (and
even $X \cap N_\lambda \backslash Z \subseteq Y^+)$; by requirements (c)
(and (a)) of $\circledast$ clearly $|Y^+| \le \lambda$, by requirement (a) of
$\circledast$ clearly $Y^+ \subseteq X \backslash Z$ and by rquirement
(b) of $\circledast$, even $Y^+ = X \cap N_\lambda \backslash Z$.

By requirements (a) + (d) of $\circledast$, clearly $h$ is a function from
$Y^+$ to $\mu$.
So in $(*)_\lambda$ for $Y,Z$ demands (a),(b),(c) on $Y^+,h$ 
are satisfied so it suffice to prove
demand $(d)$ there. So suppose $A \in {\mathcal P}^*_{\mathbf x}, \kappa^\ast
\le \theta < \mu$ and moreover, $\theta$ witness that the triple
$(\kappa,\kappa^*,\mu)$ is $\lambda$-admissible, $|A \cap Y^+| \ge \theta$
and $|A \cap Z| < \mu$ and $\beta < \mu$; we should prove 
``for every $\beta < \mu,h \restriction (h^2_B)^{-1}(\{\beta\}) 
\cap Y^+)$ coincides with $h^1_B$ on a set of cardinality $\mu$ for 
some $B\in {\mathcal Q}^*_{{\bf x},A}$".
So $|A \cap N_\lambda| \ge \theta$.  Choose a pair
$(\delta^\ast,\theta^\ast)$ such that:
\mn
\begin{enumerate}
\item[$\otimes$]
\begin{enumerate}
\item[(i)]  $\delta^\ast \le \lambda$,
\sn
\item[(ii)] $\theta^*$ witnesses that $(\kappa,\kappa^*,\mu)$ is
 $(\mu + |\delta^*|)$-admissible hence $\kappa^\ast \le 
\theta^\ast < \mu$,
\sn
\item[(iii)] $|A \cap N_{\delta^\ast}|\ge \mu$ or $\delta^\ast = \lambda$,
\sn
\item[(iv)]  under (i) + (ii) + (iii), $\delta^\ast$ is minimal.
\end{enumerate}
\end{enumerate}
\mn
This pair is well defined as  $(\lambda,\theta)$  satisfies requirement
(i) + (ii) + (iii).
\medskip

\noindent{\em Subcase B1:}\quad  $\delta^\ast$ is zero.
\medskip

So  $|Y^+_0 \cap  A| \ge \theta^\ast \ge \kappa^\ast$ hence by the
choice of $h_0$, i.e. clause (e) of $\circledast$, recalling $A \in
{\mathcal P}^*_{\mathbf x}$ we are done.
\medskip

\noindent{\em Subcase B2:}\quad  $\delta^\ast = i + 1$.
\medskip

So $\delta^*<\lambda$; clearly the pair $(i,\theta^*)$ 
standing for $(\delta^*,\theta^*)$
satisfies clauses (i)+(ii) of $\otimes$ so it cannot satisfies clause
(iii) there as then $(\delta^*,\theta)$ fail clause (iv). This means that
$|A \cap  N_i| < \mu$, but $\bigcup \{Y^+_j:j<i\}\subseteq N_i$ hence
$|A \cap\bigcup\limits_{j < i} Y^+_j| < \mu$, but also $|A \cap Z| <
\mu$ hence  $|A \cap (Z \cup\bigcup\limits_{j < i}Y^+_j)| < \mu$.
Clearly $\theta^*$ witness $(\kappa,\kappa^*,\mu)$ is $(\mu + |i|)$
admissible holds (as $\mu + |i| = \mu + |i+1| =
\mu + |\delta^\ast|)$,  so if  $|A \cap
Y^+_i| \ge \theta^\ast$ we are done by the choice of  $h_i$, i.e. by clause
(e) of $\circledast$;  if not, then $|A \cap (Z \cup \bigcup\limits_{j <
  i+1}Y^+_j)| < \mu$ and so necessarily 
$A \cap  Y^+_{i+1} \supseteq A \cap N_{i+1} \setminus \cup \{Y^+_j:j<i+1\} = A
\cap N_{\delta ^\ast} \setminus
\cup \{Y^+_j:j<i+1\}$ has cardinality  $\ge \theta^\ast$ (and ``$\theta^*$
witness $(\kappa,\kappa^*,\mu)$ is $\mu + |i+1| = |Y^+_{i+1}|$--admissible"
holds) so we are done by the choice of  $h_{i+1}$.
\medskip

\noindent{\em Subcase B3:}\quad  $\delta^*$ is a limit ordinal below
$\lambda$.
\medskip

So for some $i < \delta^\ast$,  $|A \cap N_i| \ge \theta^\ast$. [Why? As
$\theta^\ast < \mu \le |A \cap N_{\delta^*}|$].  Now in $N_{i+1}$ there is a
maximal family  ${\mathcal Q} 
\subseteq [X \cap N_i]^{\theta^\ast}$ satisfying $[B_1
\ne B_2 \in {\mathcal Q} \Rightarrow  |B_1 \cap B_2| < \kappa^\ast]$
hence by clause
(ii) of $\otimes$ and clause (c) of Definition \ref{ref.45}(2)
we have $|{\mathcal Q}| \le \mu + |\delta^\ast|$.  Choosing the
$<^*_\chi$-first such ${\mathcal Q}$, 
clearly ${\mathcal Q} \in N_{i+1}$ so recalling ${\mathcal Q} 
\in N_{i+1} \subseteq N_{\delta^*}$ we have $\cQ 
\subseteq N_{\delta^*}$.  By the choice of 
${\mathcal Q}$, necessarily there is $B \in {\mathcal Q}$ such
that $|B \cap  A| \ge \kappa^\ast$ (if $A \notin \cQ$ by the
maximality of $\cQ$ and if $A \in \cQ$ one can choose $B = A$), 
but as $B \in {\mathcal Q}$ clearly $B \in
N_{\delta^*}$ and $|B| = \theta^\ast < \mu  = \mu^\kappa$
hence  $\left[ B^\prime \in  [B \cap A]^\kappa \Rightarrow  B \cap A \in
N_{\delta^\ast} \right] $.  As $A \in {\mathcal P}^\ast_{\mathbf x}$
and $|B \cap A| \ge \kappa^*$ 
there is $B^\prime \in [B \cap A]^\kappa$ satisfying
$c\ell(B') \subseteq A$, $|c \ell(B')|  = \mu$.
Clearly $c \ell(B') \in  N_{\delta^\ast}$ hence for some
$j \in (i,\delta^\ast)$, $c \ell(B') \in N_j$ hence $c \ell(B') \subseteq
X \cap N_j$.  So  $|A \cap N_j| \ge \mu$.  By assumption for some
$\theta^\prime \in [\kappa^\ast,\mu)$ the triple $(\kappa,\kappa^*,\mu)$ is
$(\mu + |j|)$-admissible, see Definition \ref{ref.45}, so the pair
$(j,\theta^\prime)$  contradicts the choice of $(\delta^\ast,\theta^\ast)$.
\medskip

\noindent
{\em Subcase B4:} \quad  $\delta^*=\lambda$.
\medskip

As  $\lambda \in N_0$,  there is a maximal family  ${\mathcal Q} \subseteq
[\lambda]^{\theta^\ast}$ satisfying

\[[B_1 \ne B_2 \in {\mathcal Q}\quad 
\Rightarrow\quad |B_1 \cap B_2|< \kappa^\ast]\]

which belongs to  $N_0$.  By the assumption $\otimes(ii)$ on
$\theta^*$ and clause (c) of Definition \ref{ref.45}(2) we
know that $|{\mathcal Q}| \le \lambda$, but $\lambda+1\subseteq N_\lambda$ 
hence ${\mathcal Q} \subseteq N_\lambda$ hence $(\forall B \in {\mathcal Q})
(\exists i < \lambda)(B \in N_i)$.  We define by induction on $j \le
\lambda$,  a one-to-one function  $g_j$ from  $N_j \cap X \backslash
Z$  onto an initial segment of  $\lambda$ increasing continuous with $j$,
$g_j$ the $<^\ast_\chi$-first such function.  So clearly $g_j \in N_{j+1}$
and let ${\mathcal Q}^\prime = \{g^{-1}_\lambda (B):B \in  {\mathcal Q}\}$,
(i.e. $\{ \{ g^{-1}_\lambda (x):x \in B \}:B \in 
{\mathcal Q}\}$). Clearly for any $B \in {\mathcal Q}$, there is $i <
\lambda$ such that $B \in N_i \cap {\mathcal Q}$, let $\mathbf i(B)$
be the first such 
$i$, so $B \subseteq {\rm  Dom}(g^{-1}_{\mathbf i(B)})$ and so 
$g^{-1}_{\mathbf i(B)}(B) \in N_{i+1}$ and $g_\lambda$ is necessarily 
a one-to-one function from
$N_\lambda \cap X\backslash Z$  onto  $\lambda$. Recall that
$A \cap Y^+=A \cap (X \cap N_\lambda) \backslash Z$ has cardinality
$\ge \theta^*$.  Hence for some  $B \in
{\mathcal Q}^\prime$,  $|B \cap  A| \ge \kappa^\ast$, so as in subcase B3, for
some  $B^\prime \in N_\lambda$,  $B^\prime \subseteq B \cap A$,  $|B^\prime|
=  \kappa$, $c\ell(B') \subseteq A$, $|c\ell(B')| = \mu$.  Clearly $B\in
N_{\mathbf i(B)+1}$ hence $[B]^{\le \kappa} \in N_{{\mathbf i}(B)+1}$ but
its cardinality is $\le \mu$ hence $[B]^{\le \kappa}\subseteq N_{{\mathbf
    i}(B)+1}$, so $B' \in N_{\mathbf i(B)+1}$ and so $c\ell(B')
\subseteq N_{{\mathbf i}(B)+1}$. But $|A \cap Z| < \mu$ so 
by the last two sentences $|A \cap
Y^+_{{\mathbf i}(B)+ 1}| = \mu$ and by assumption $\otimes$(ii), some
$\theta$ is a witness to $(\kappa,\kappa^\ast,\mu)$ being $(\mu
+|i|)$-admissible (stipulating $i =  \bold i(B)+1$), 
contradicting the choice of $(\delta^*,\theta^*)$
(i.e. minimality of $\delta^\ast$).  \hfill$\square_{\ref{ref.49}}$

\begin{discussion}
\label{3.7}
1) If we would like to include the case 
$\mu = 2^{\aleph_0} = \aleph_2$, $\kappa = \aleph_0,\kappa^* =\aleph_1$
we should consider a revised framework.  We have a family
${\mathfrak I}$ of ideals on cardinals $\theta$ from $[\kappa^*,\mu)$ 
which are $\kappa$-based
(i.e. if $A \in I^+$, $I \in {\mathfrak I}$ (similar to \cite{HJSh:249})
then $(\exists B \in [A]^\kappa)(B \in I^+$)) and in Definition
\ref{ref.48}(3) hence in the proof of \ref{ref.49} replace ${\mathcal
P}^*_{\bold x}$ by 

\[\begin{array}{ll}
{\mathcal P}^\ast={\mathcal P}^\ast_{\mathfrak I}=:\Big\{A\subseteq X:
&|A| = \mu\mbox{ and for every pairwise distinct } x_\alpha \in A 
\mbox{ for} \\ 
&\alpha<\theta \mbox{ the set }\{u\subseteq \theta:
|c\ell\{x_\alpha:\alpha\in u\}|< \mu\}\\
&\mbox{is included in some } I \in {\mathfrak I} \Big\}.
\end{array}\]
\medskip

\noindent
and in Definition \ref{ref.45}(1),(2) we replace the triple
$(\kappa,\kappa^*,\mu)$ by the quadruple $(\kappa,\kappa^*,\mu,\gI)$
and clause (c) of \ref{ref.45}(2) by
\mn
\begin{enumerate}
\item[$(c)'_\lambda$]   $\lambda \ge \mu$ and: $|\cF| \le \lambda$ whenever
\[\cF \subseteq \{ (\theta,I,f):I \in {\mathfrak I},\
\theta = {\rm  Dom}(I),\ f:\theta\longrightarrow \lambda
\mbox{ is one to one}\},\]
and if $(\theta_\ell,I_\ell,f_\ell) \in \cF$ for $\ell = 1,2$ are distinct
then
\[\{ \alpha < \theta_2:f_2(\alpha) \in {\rm  Rang }(f_1)\} \in I_2.\]
\end{enumerate}
\mn
Note that the present ${\mathcal P}^\ast$ fits for 
 repeating the proof of \ref{ref.49}.

\noindent
2) \underline{Discussion of the Consistency of {\sc NO}}:

There are some restrictions on such theorems.  Suppose
\mn
\begin{enumerate}
\item[$(*)$] GCH and there is a stationary
$S \subseteq \{\delta <\aleph_{\omega +1}:{\rm cf}(\delta) = \aleph_1\}$ and
$\langle A_\delta:\delta \in S \rangle$ such that:
\begin{itemize}
\item $A_\delta \subseteq \delta = \sup A_\delta$,
\sn
\item ${\rm otp}(A_\delta) = \omega_1$ and
\sn
\item $\delta_1 \ne \delta_2\quad \Rightarrow\quad |A_{\delta_1} \cap
  A_{\delta_2}| < \aleph_0$.
\end{itemize}
\end{enumerate}
\mn
(This statement is consistent by \cite[4.6,p.384]{HJSh:249} which continues
\cite{Sh:108} see more in \cite{Sh:668}.)

Now on $\aleph_{\omega +1}$ we define a closure operation:

\[
\alpha \in c \ell(u)\quad \Leftrightarrow\quad (\exists \delta \in
S)[\alpha \in A_\delta \mbox{ and } (u \cap A_\delta) \ge \aleph_0].
\]

This certainly satisfies the demands in Definition \ref{ref.48} with $\kappa =
\kappa^* = \aleph_0,\mu = \aleph_1$ except the pcf assumptions,
i.e. clause (c) of Definition \ref{ref.45}(2).  
However, this is not a case of our theorem.

\noindent
3) We may consider in the proof of \ref{ref.49} strengthening clause (e) of
$\circledast$ by weakening clause $(e)(\delta)$ of $\circledast$ by
fixing the ordinal $\beta$ and demanding only $(A \backslash
\bigcup\limits_{j<i} Y^+_j \backslash Z) \cap
(h^z_\beta)^{-1}(\{\beta\})$ has cardinality $\mu$.  But we do not
seem to gain anything.
\end{discussion}


\end{document}